\documentclass[12pt]{amsart}
\usepackage{amsthm,amsbsy,amsfonts,amssymb,amsmath,amscd}

\newcommand{\Q}{\mathbb Q}
\newcommand{\Z}{\mathbb Z}
\newcommand{\R}{\mathbb R}
\newcommand{\C}{\mathbb C}
\newcommand{\G}{\mathcal G}

\newcommand{\nN}{\mathcal N}
\newcommand{\A}{\mathbb A}

\newcommand{\s}{\text sym}

\makeatletter
\def\l@section{\@tocline{1}{4pt}{1pc}{}{}}
\def\l@subsection{\@tocline{2}{0pt}{2pc}{5pc}{}}
\makeatother

\makeatother 
\title{Icosahedral fibres of the symmetric cube and algebraicity}

\author{Dinakar Ramakrishnan}\thanks{Partly supported by the NSF grant DMS-0701089}

\begin{document}
\subjclass[2000]{11F70; 11F80; 22E55}
\maketitle

\medskip

\begin{flushright}
{\it To Freydoon Shahidi \\ On the occasion of his sixtieth birthday \\
Tavalodet Mubarak!}
\end{flushright}

\bigskip

\section*{\bf Abstract}

For any number field $F$, call a cusp form
$\pi=\pi_\infty\otimes\pi_f$ on GL$(2)/F$
{\it special icosahedral}, or just
{\it s-icosahedral} for short, if $\pi$ is not solvable
polyhedral, {\it and} for a suitable ``conjugate'' cusp form $\pi'$
on GL$(2)/F$, sym$^3(\pi)$
is isomorphic to sym$^3(\pi')$, and the symmetric fifth power
$L$-series of $\pi$ equals the Rankin-Selberg $L$-function $L(s,
{\rm sym}^2(\pi')\times \pi)$ (up to a finite number of Euler
factors). Then the point of this Note is to obtain the
following result:

\it Let $\pi$ be s-icosahedral (of trivial central
character). Then $\pi_f$ is algebraic without local components of
Steinberg type, $\pi_\infty$ is of Galois type, and $\pi_v$ is
tempered everywhere.
Moreover, if $\pi'$ is also of trivial central character,
it is s-icosahedral, and the field of
rationality $\Q(\pi_f)$ (of $\pi_f$) is $K:=\Q[\sqrt{5}]$, with
$\pi'_f$ being the Galois conjugate of $\pi_f$ under
the non-trivial automorphism of $K$.\rm \,

There is an analogue in the case of non-trivial
central character $\omega$, with the conclusion that $\pi_f$
is algebraic when
$\omega$ is, and when $\omega$ has finite order, $\Q(\pi_f)$ is contained
in a cyclotomic field.

\medskip

\section*
{\bf Introduction}

\bigskip

Let us begin with some motivation and consider a continuous
irreducible representation $\rho$ with trivial determinant of the
absolute Galois group $\G_F$ (of a number field $F$) into GL$(2,\C)$
which is {\it icosahedral}, i.e., whose image in PGL$(2,\C)$ is the
alternating group $A_5$. Then it is well known that $\rho$ is
rational over $K=\Q[\sqrt{5}]$, but not equivalent to its Galois
conjugate $\rho^\prime$ under the non-trivial automorphism of $K$.
Moreover, one has (see \cite{K2}, \cite{Wang}, for example):  (i)
$\s^3(\rho) \simeq \s^3(\rho^\prime)$, and (ii) $\s^5(\rho) \simeq
\s^2(\rho^\prime)\otimes\rho$. These two are the signature
properties of such representations, and they lend
themselves to natural automorphic analogues.

Let $\pi$ be a cuspidal automorphic
representation of GL$(2,\A_F)$ with central character $\omega$. For every $m\geq
1$ one has its {\it symmetric $m$-th power $L$-function}
$L(s,\pi;{\rm sym}^m)$, which is an Euler product over the places
$v$ of $F$, with the $v$-factors (for finite unramified $v$ of norm
$q_v$) being given by
$$
L_v(s,\pi;{\rm sym}^m) \, = \, \prod_{j=0}^m
(1-\alpha_v^j\beta_v^{m-j}{q_v}^{-s})^{-1},
$$
where the unordered pair $\{\alpha_v, \beta_v\}$ defines the
diagonal conjugacy class in GL$_2(\C)$ attached to $\pi_v$. Even at
a ramified (resp. archimedean) place $v$, one has by the local
Langlands correspondence a $2$-dimensional representation $\sigma_v$
of the extended Weil group $W_{F_v}\times{\rm SL}(2,\C)$ (resp. of
the Weil group $W_{F_v}$), and the $v$-factor of the symmetric
$m$-th power $L$-function is associated to sym$^m(\sigma_v)$. A special case of the
{\it principle of functoriality} of Langlands asserts that there is, for each $m$,
an (isobaric) automorphic representation ${\rm sym}^m(\pi)$ of
GL$(m+1,\A)$ whose standard (degree $m+1$) $L$-function $L(s, {\rm
sym}^m(\pi))$ agrees, at least at the primes not dividing $\nN$,
with $L(s,\pi;{\rm sym}^m)$. This was known to be true for $m=2$ long ago by the
work of Gelbart and Jacquet (\cite{GeJ}), and more recently for $m= 3, 4$ by the
deep works of Shahidi and Kim (\cite{KSh1, K, KSh2}). Write $A^{2j}(\pi)$ for
sym$^{2j}(\pi)\otimes\omega^{-j}$ when the latter is automorphic;
it is also customary to denote $A^2(\pi)$ by Ad$(\pi)$.
We will say that $\pi$ is {\it solvable polyhedral} if
sym$^m(\pi)$ is Eisensteinian for some $m \leq 4$.

\medskip

Suppose $\pi'$ is another cusp form on GL$(2)/F$, say of the same
central character as $\pi$, such that sym$^2(\pi) \simeq {\rm sym}^2(\pi')$. Then one
knows (\cite{Ra2}) that $\pi$ and $\pi'$ must be abelian twists of each other. One could
ask if the same conclusion holds if sym$^3(\pi)$ is isomorphic to sym$^3(\pi')$. The
answer is in the negative in that case, and a counterexample would be furnished by a
$\pi$ associated to a $2$-dimensional icosahedral Galois representation $\rho$
(of trivial determinant), i.e.,
with its image in PGL$(2,\C)$ being isomorphic to the
alternating group $A_5$. Indeed, as remarked above, $\rho$ would be defined
over $\Q[\sqrt{5}]$ and $\pi'$ would be associated to the Galois conjugate $\rho'$ of
$\rho$ (under $a+b\sqrt{5}\mapsto a-b\sqrt{5}$). However, the {\it even} Galois representations
are not (at all) known to be modular. Nevertheless, cup forms of such icosahedral type
are of great interest to study, for themselves and for understanding the fibres of the
symmetric cube transfer. What we do is give a definition of an s-icosahedral cusp
form which does not depend on any conjecture, and which is robust enough to
furnish consequences which one would usually know only when there is an associated
Galois representation.

\medskip

Let us call a cusp form $\pi$ on GL$(2)/F$ (of central character $\omega$) {\it s-icosahedral} if we have
\begin{enumerate}
\item[(sI-1)] $\pi$ is not solvable polyhedral; \,
\item[(sI-2)] sym$^3(\pi) \simeq {\rm sym}^3(\pi')$,
for a cusp form $\pi'$ on GL$(2)/F$; \, and
\item[(sI-3)] the following identity of $L$-functions holds (outside a finite set $S$ of
places of $F$ containing the archimedean and ramified places):
$$
L^S(s, \pi; {\rm sym}^5) \, = \, L^S(s, {\rm Ad}(\pi')\times \pi\otimes\omega^2),
$$
with $\pi'$ as in (sI-2).
\end{enumerate}

\medskip

Observe that if $(\pi, \pi')$ are associated as above, then $(\pi, \pi'\otimes\nu)$ will also be associated
if $\nu$ is a cubic character, but this ambiguity can be eliminated by requiring that the ratio of the central
characters $\omega, \omega'$, of $\pi, \pi'$ respectively, is not cubic. Note also that the property of being
s-icosahedral is invariant under twisting $\pi \mapsto \pi \otimes \chi$ by
idele class characters $\chi$, with $\pi'\mapsto \pi'\otimes\chi$, in particular under taking contragredients
$\pi\mapsto \pi^\vee=\pi\otimes\omega^{-1}$.

\medskip

Given cusp forms $\pi_1, \pi_2$ on GL$(n_1)/F$, GL$(n_2)/F$ respectively, one has good
analytic properties (\cite{Sh1, Sh3, JPSS, JS1}) of the
associated Rankin-Selberg $L$-function $L(s, \pi_1 \times \pi_2)$. When
$(n_1,n_2)=(2,2)$ (\cite{Ra2}) and $(n_1, n_2)=(3,2)$ (\cite{KSh1}), one also knows
the existence of an isobaric automorphic form $\pi_1\boxtimes\pi_2$
on GL$(n_1n_2)$.
Thus the hypothesis (sI-3) implies that sym$^5(\pi)$ is
modular; in fact, ${\rm Ad}(\pi')\boxtimes \pi\otimes\omega^2$ represents it at every
finite place (and at infinity), as seen by the standard stability results due to Shahidi and others
(see section 2 below).

\medskip

What we want to do in this Note is to look at the
situation where one does not know of the existence of a corresponding Galois representation $\rho$, like when
$\pi$ is a suitable Maass form on GL$(2)/\Q$, and see what arithmetic properties one can still deduce
{\it unconditionally}. However it is a {\it Catch $22$} situation because it is not easy, either by the trace
formula or by other means, to construct such Maass forms of Laplacian eigenvalue $1/4$ except by starting with
a Galois representation. Nevertheless, here is what we prove:

\medskip

\noindent{\bf Theorem A} \, \it Let $F$ be a number field and
$\pi=\pi_\infty\otimes\pi_f$ a cuspidal automorphic representation
of GL$_2(\A_F)$ with algebraic central character $\omega$.
Suppose $\pi$ is s-icosahedral, i.e., satisfies (sI-1), (sI-2) and
(sI-3) relative to a cusp form $\pi'$ on GL$(2)/F$ of central character $\omega'$. Then we have the
following:
\begin{enumerate}
\item[(a)] $\pi_f$ is algebraic with no local components of Steinberg type
\item[(b)] If $\omega$ is of finite order, $\pi_f$ is rational over a cyclotomic
extension of $\Q$, and $\pi_\infty$ is of Galois type
\item[(c)] $\pi$ is tempered, i.e., satisfies the Ramanujan hypothesis.
\item[(d)] If $\omega=\omega'=1$, we have:
\begin{enumerate}
\item[(1)]$\pi_f$ is rational over $\Q[\sqrt{5}]$, but not over $\Q$;
\item[(ii)]$\pi'$ is unique, with $\pi'_f$ being
the Galois conjugate of $\pi_f$ under the non-trivial
automorphism $\tau$ of $\Q[\sqrt{5}]$; \, and
\item[(iii)]$\pi'$ is also s-icosahedral and
of Galois type.
\end{enumerate}
\end{enumerate}
\rm

\medskip

Concerning (b), let us recall that $\pi_\infty$ is said to be of {\it Galois type} if at
every archimedean place $v$,
the associated $2$-dimensional representation
$\sigma_v$ of the Weil group $W_{F_v}$ is trivial upon restriction to $\C^\ast$.

One writes $\Q(\pi_f)$ for the field of rationality of $\pi_f$ (\cite{C}, \cite{Wa}). We will use
the algebraic parametrization at any $v$ which is equivariant for the action of Aut $\C$ on the two
sides of the local Langlands correspondence (\cite{He2}. In particular, for $v$ unramified, it corresponds to
the {\it Tate paramaterization} of \cite{C}, which was introduced so as not to worry about
spurious square-roots of $q_v=Nv$; the normalization is unitary in \cite{Wa}.

\medskip

In \cite{Ra7}, we introduced a notion of quasi-icosahedrality, based
on a conditional assumption that certain symmetric powers of $\pi$ are automorphic.
More precisely, an irreducible cuspidal automorphic representation $\pi$ of GL$(2,\A_F)$
is called {\it quasi-icosahedral} iff we have
\begin{enumerate}
\item[(i)]sym$^m(\pi)$ is automorphic for every $m \leq 6$;
\item[(ii)]sym$^m(\pi)$ is cuspidal for every $m \leq 4$;
and
\item[(iii)]sym$^6(\pi)$ is not cuspidal.
\end{enumerate}
The key result of \cite{Ra7} (see part (b) of Theorem A$^\prime$ of section 2) is that, for every such quasi-icosahedral $\pi$,
there exists another
cusp form $\pi'$ of GL$(2)/F$ such that the symmetric fifth power
of such a quasi-icosahedral cusp form $\pi$ is necessarily a character twist of the functorial
product ${\rm Ad}(\pi')\boxtimes\pi$. From this we obtain

\medskip

\noindent{\bf Proposition B} \, \it Let $\pi$ be an s-icosahedral cusp form on
GL$(2)/F$ of central character $\omega$. Assume that sym$^6(\pi)$ is automorphic. Then $\pi$ is quasi-icosahedral.

\rm

\medskip

It will be left to the astute reader to figure out the conditions on a quasi-icosahedral $\pi$ which will make it s-icosahedral.
We do not pursue this any further here because we want to stay in the unconditional realm here.

\medskip

Sometimes in Mathematics, when one makes a right prediction and puts
it in a workable framework, the proof is not hard to find; Theorem A here is an instance of that
and the proof mostly requires just
book-keeping. Nevertheless, we hope that the conclusion is of some
interest. It should be mentioned that given the definition, it is not surprising that an
s-icosahedral $\pi$ is algebraic,
but what is nice is that when $\omega=1$, $\pi$ is even rational over $\Q(\sqrt{5})$.

\medskip

Finally, one can easily construct non-singleton fibres of the symmetric cube transfer by taking, 
for any cusp form $\pi$ on GL$(2)$, the collection $\{\pi\otimes\chi \, \vert \, \chi^3=1\}$. 
(By \cite{KSh2}, sym$^3(\pi)$ will be cuspidal iff $\pi$ is not dihedral or tetrahedral.) These 
fibres are not terribly interesting, however, especially compared to the icosahedral ones.

\medskip

We thank the referee for helpful comments, and the NSF for continued support through the grant DMS-0701089.

\medskip

This article is dedicated to Freydoon Shahidi, from whom I have
learnt a lot about automorphic $L$-functions and the Langlands-Shahidi
method, and who has been a
longtime friend and in effect a {\it baradar e bozorgtar}.

\vskip 0.2in

\section{\bf The symmetric cube constraint}

\medskip

Let $F$ be a number field with adele ring $\A_F$. For each place $v$, denote by
$F_v$ the corresponding local completion of $F$,
and for $v$ finite, by $\mathfrak O_v$ the ring of integers of $F_v$ with uniformizer
$\varpi_v$ of norm $q_v$. Throughout this article, $\omega$ (resp. $\omega'$) denote the central character
of a cusp form $\pi$ (resp. $\pi'$) on GL$(2)/F$.

We will use, without mention, the notations and conventions of
\cite{Ra7}, especially section 1 therein.

\medskip

The symmetric cube condition (sI-2) (see Introduction) imposes strong constraints on how the local
components $\pi_v$ and $\pi'_v$ are related at all $v$. We spell these out at the
unramified $v$ in the second part of the Lemma
below, and show that, not surprisingly, when the central character $\omega$ is algebraic, $\pi_v$
and $\pi'_v$ are algebraically related. The first part show the effect of $\pi$ satisfying both (sI-1)
and (sI-2) on its ``mirror'' $\pi'$.

\bigskip

\noindent{\bf Lemma 1.1} \, \it Let $\pi$ be a cusp form on
GL$(2)/F$ which is s-icosahedral
relative to $\pi'$. Then
\begin{enumerate}
\item[(A)]$\pi'$ satisfies (sI-j) for $j=1,2$, and moreover,
$$
(\omega')^3 \, = \, \omega^3.
$$
In particular, $\omega'$ is algebraic iff $\omega$ is.
\item[(B)]Let $v$ be any finite place where $\pi, \pi'$ are
unramified, and denote by $\{a,b\}$ with $b=wa^{-1}$ (resp. $\{c,d\}$
with $d=w'c^{-1}$) the unordered pair of complex numbers associated to $\pi_v$ (resp.
$\pi'_v$), where $w=\omega_v(\varpi), w'=\omega'_v(\varpi)$. Write
$$
w=zw', \, \, \, {\rm with} \, \, \, z^3=1.
$$
Then one of the following cases occurs (up to interchanging \newline$a$ and $b$):
\begin{enumerate}
\item[(1)]$\{c, d\} = \{za, zb\}$; Ad$(\pi_v) \simeq {\rm Ad}(\pi'_v)$
\item[(2)]$\{c, d\}=\{\mu za, \mu^{-1} zb\}$, $\mu^4=1$, $a^2=\mu w$; Ad$(\pi_v) \simeq {\rm Ad}(\pi'_v)$
\item[(3)]$\{c, d\} = \{\zeta za, \zeta^{-1} zb\}$, $\zeta^5=1$, $a^2=\zeta w$
\end{enumerate}
In particular, when $\omega$ is algebraic, $\pi_v$ and $\pi'_v$ are both
algebraic in cases $(2)$ and $(3)$. Moreover, if $\omega$ has finite order, $\Q(\pi_v)$
is contained, in these cases, in a finite abelian extension $K$ of $\Q$, independent of $v$,
containing $\Q(\omega)$.
\end{enumerate}
\rm

\medskip

\noindent{\bf Remark 1.2} \, In section 3 we will show, by also appealing to (sI-3), which is not used in Lemma 1.1, that (i) $\pi_v$ is algebraic
in {\it all} cases, and (ii) $a^4=1$ in case (2) of part (B) above, if $w=1$. It is perhaps useful to note that if (sI-3) is not satisfied by $\pi$, one may take $F=\Q$ and $\pi'$ to be a cubic twist of $\pi$, and such a $\pi$ is expected to be
transcendental if it is generated, for example, by a Maass form $\varphi$ of weight $0$ relative to a congruence subgroup of SL$(2,\Z)$ acting on the upper half plane with Laplacian eigenvalue $\lambda > 1/4$.

\medskip

{\it Proof of Lemma 1.1} \, (A) \, The fact that $\pi$ and $\pi'$ have isomorphic symmetric cubes implies immediately that ${\omega'}^6=\omega^6$, which is not sufficient for us. We will first show that $\pi'$ is not solvable polyhedral. First, since $\pi$ is not solvable polyhedral, sym$^j(\pi)$ is cuspidal for $j\leq 4$ (cf. \cite{KSh2}). By (sI-2), sym$^3(\pi')$ is also cuspidal. Suppose
sym$^4(\pi')$ is not cuspidal. Then, by the criterion of Kim and Shahidi (\cite{KSh2}), sym$^3(\pi')$ must be monomial, which forces it to admit a non-trivial self-twist by a quadratic character $\nu$, say. Then by (sI-2), sym$^3(\pi)$ also admits a self-twist by $\nu$, implying that sym$^4(\pi)$ is not cuspidal, contradicting (sI-1). Hence $\pi'$ satisfies both (sI-j) for $j \leq 2$.

Next we appeal to Kim's theorem (\cite{K}) giving the automorphy in GL$(6)/F$ of the exterior square $\Lambda^2(\Pi)$ of any cusp form $\Pi$ on GL$(4)/F$. Applying this with $\Pi={\rm sym}^3(\eta)$ for a cusp form $\eta$ on GL$(2)/F$ with central character $\nu$, we obtain
the (well known) isobaric decomposition
$$
\Lambda^2({\rm sym}^3(\eta)) \, = \, \left({\rm sym}^4(\eta)\otimes \nu\right) \boxplus \nu^3.\leqno(1.3)
$$
To prove this, we first note that by \cite{K}, the $L$-functions agree at almost all places, and then appeal to the strong multiplicity one theorem for global isobaric representations, due to Jacquet and Shalika (\cite{JS1}).
Now applying this to $\eta=\pi, \pi'$, we get by (sI-2), the following equivalence of isobaric sums:
$$
\left({\rm sym}^4(\pi)\otimes \omega\right) \boxplus \omega^3 \, = \, \left({\rm sym}^4(\pi')\otimes \omega'\right) \boxplus {\omega'}^3.\leqno(1.4)
$$
Since ${\rm sym}^4(\pi)$ and ${\rm sym}^4(\pi')$ are both cuspidal, we are forced to have
$$
(\omega')^3 \, = \, \omega^3,
$$
as claimed.

\medskip

(B) \, Preserving the notations used in part (B) of the Lemma, and noting that
$$
b=wa^{-1}, \, d=w'c^{-1},
$$
we get by (sI-2),
the equality of sets:
$$
\{a^3, a^2b, ab^2, b^3\} \, = \, \{c^3, c^2d, cd^2, d^3\}.\leqno(1.5)
$$
Clearly,
$$
a^2b=wa, \, ab^2=w^2a^{-1}, b^3=w^3a^{-3} \, \, \, {\rm and} \, \, \, c^2d=w'c, \, cd^2={w'}^2c^{-1}, \, d^3={w'}^3c^{-3}.\leqno(1.6)
$$

Note that, since $w'=z^2w$, we have
$$
c=\alpha a \, \implies \, d=w'\alpha^{-1}a^{-1} = z^2\alpha^{-1}b,\leqno(1.7)
$$
and
$$
c=\beta b \, \implies \, d=w'\beta^{-1}b^{-1} = z^2\beta^{-1} a.\leqno(1.8)
$$

\bigskip

{\it A priori}, $a^3$ has four possibilities to satisfy (1.5). However, by interchanging $c$ and $d$, we are reduced to
considering only {\sl two main cases}:

\newpage

\noindent{{\it Case {\bf I}}: \, \, ${\bf a^3=c^3}$} \,

\medskip

In this case, since ${\omega'}^3=\omega^3$ and $w'=z^2w$, (1.5) and (1.6) yield, after dividing by $w$,
$$
\{a, wa^{-1}\} \, = \, \{-z^2c, zwc^{-1}\}.\leqno(1.9)
$$
If $a=-z^2c$, then $c=za$ and $d=zb$ (by (1.7)), putting us in the situation (1). Since we have
$$
\{c,d\} \, = \, \{a,b\}\cdot z,
$$
$\pi_v$ is (at most) a cubic twist of $\pi_v'$; it follows that Ad$(\pi_v)$ and Ad$(\pi'_v)$ are isomorphic.

So we may assume $a=zwc^{-1}$. Then we get (using (1.8))
$$
c = zwa^{-1}= zb, \, \, \, {\rm and} \, \, \, d=za,\leqno(1.10)
$$
again landing us in (1).

\bigskip

\noindent{{\it Case {\bf II}}: \, \, ${\bf a^3=w'c}$}

\medskip

In this case, $c={w'}^{-1}a^{3}$, and we obtain from (1.5) and (1.6),
$$
\{wa, w^2a^{-1}\} \, = \, \{w^{-3}a^{9}, {w}^{6}a^{-9}\}.\leqno(1.11)
$$
If $wa=w^{6}a^{-9}$, then $a^{10}=w^{5}$, so that
$$
a^2 \, = \, \zeta w, \, \, \, {\rm with} \, \, \, \zeta^5=1.\leqno(1.12)
$$
Hence $a^3=\zeta wa$, and so
$$
c={w'}^{-1}\zeta w a \, = \, \zeta z a.\leqno(1.13)
$$
It follows that
$$
d = w'c^{-1} = (z^2w)\zeta^{-1} z^{-1}a^{-1} \, = \, \zeta zb.\leqno(1.14)
$$
Thus we are in situation (3).

\medskip

It is left to consider when $wa=w^{-3}a^9$. In this case, $a^8 = w^4$, so that
$$
a^2 \, = \, \mu w, \, \, \, {\rm with} \, \, \, \mu^4=1.\leqno(1.15)
$$
Hence $a^3 = \mu wa$, and arguing as above, we deduce
$$
c \, = \, \mu za, \, \, \, {\rm and} \, \, \, d \, = \, \mu^{-1} zb.\leqno(1.16)
$$
This puts us in situation (2).

We still need to show that Ad$(\pi_v) \simeq {\rm Ad}(\pi'_v)$ in this case as well. To this end, note that since $a^2=\mu w$ and $c=\mu za$,
$$
c \, = \, (a^2w^{-1})za \, = \, w^{-1}za^3.\leqno(1.17)
$$
On the other hand, raising $c=\mu za$ to the fourth power, we get, since $\mu^4=1$,
$$
c^4 \, = \, za^4 \, = \, za(a^3) \, = \, za(wz^{-1}c),\leqno(1.18)
$$
which in turn yields
$$
c^3 \, = \, wza.\leqno(1.19)
$$
Furthermore,
$$
c^2 \, = \, w^{-2}z^2a^6 \, = \, w^{-2}z(za^4)a^2 \, = \, w^{-2}zc^4a^2,
$$
implying
$$
c^{-2} \, = \, w^{-2}za^2, \, \, \, {\rm and} \, \, \, c^2 \, = \, w^2z^2a^{-2}.\leqno(1.20)
$$
The unramified representation Ad$(\pi_v)$ is described by the unordered triple
$$
\{ab^{-1}, 1, a^{-1}b\} \, = \, \{w^{-1}a^2, 1, wa^{-2}\}.\leqno(1.21)
$$
Similarly, Ad$(\pi_v')$ is given by the unordered triple
$$
\{cd^{-1}, 1, c^{-1}d\} \, = \, \{zw^{-1}c^2, 1, z^{-1}wc^{-2}\},
$$
which, by (1.20), is the same as
$$
\{wa^{-2}, 1, w^{-1}a^2\}.\leqno(1.22)
$$
The assertion follows.

\qed

\bigskip

\section{\bf The nicety of sym$^5(\pi)$ for $\pi$ s-icosahedral}

\medskip

As mentioned in the Introduction, the condition (s-I3) implies that for an s-icosahedral $\pi$ of
central character $\omega$, the automorphic
representation Ad$(\pi')\boxtimes\pi\otimes\omega^2$ of GL$_6(\A_F)$, whose existence is given by
\cite{KSh1}, represents sym$^5(\pi)$ at all places outside a finite set $S$ of places containing the
archimedean and ramified places. In fact we have the following strengthening:

\noindent{\bf Proposition 2.1} \, \it Let $\pi$ be an s-icosahedral cusp form on GL$(2)/F$ with central character $\omega$.
\begin{enumerate}
\item[(a)]At every finite place $v$,
$$
L(s, \pi_v; {\rm sym}^5) \, = \, L(s, {\rm Ad}(\pi'_v)\boxtimes \pi_v\otimes\omega_v^2)\leqno(2.2)
$$
and
$$
\varepsilon(s, \pi_v; {\rm sym}^5) \, = \, \varepsilon(s, {\rm Ad}(\pi'_v)\boxtimes \pi_v\otimes\omega_v^2)\leqno(2.3)
$$
\item[(b)]If $\Sigma_\infty$ denotes the set of archimedean places of $F$, we have
$$
L(s, \pi_\infty; {\rm sym}^5) \, = \, L(s, {\rm Ad}(\pi'_\infty)\boxtimes \pi_\infty\otimes\omega_\infty^2),\leqno(2.4)
$$
where $\pi_\infty$ (resp. $\pi'_\infty$) denotes $\otimes_{v\in \Sigma_\infty} \pi_v$
(resp. $\otimes_{v\in \Sigma_\infty} \pi'_v$).
\item[(c)]If $\pi, \pi'$ are not dihedral twists of each other, i.e., if their symmetric squares are not twist equivalent, then
sym$^5(\pi)$, defined to be Ad$(\pi')\boxtimes\pi\otimes\omega^2$, is cuspidal, and hence $L(s,\pi; {\rm sym}^5)$ is entire in this case.
\end{enumerate}
\rm

\bigskip

\noindent{\bf Remark} \, Of course we are dealing with a very special class of cusp forms $\pi$ here. For the status of results in the general case (which is much more complicated), see \cite{Sh1}, \cite{Sh4} and \cite{KSh2}.

\medskip

{\it Proof}. At any place $v$, we have the well known factorization (using Clebsch-Gordon):
$$
L(s, {\rm sym}^4(\pi_v)\times\pi_v) \, = \, L(s, \pi; {\rm sym}^5)L(s, {\rm sym}^3(\pi_v)\otimes\omega_v).\leqno(2.5)
$$
Similarly for the $\varepsilon$-factors, making use of the local Langlands conjecture for GL$(n)$, now known by the works of
Harris-Taylor \cite{HaT} and Henniart \cite{He}. So it holds for the $\gamma$-factors as well.
We claim that one has, for a sufficiently ramified character $\nu_v$ of $F_v^\ast$,
$$
L(s, \pi; {\rm sym}^5\otimes \nu_v) \, = \, 1.\leqno(2.6)
$$
Indeed, by the standard stability results for the Rankin-Selberg $L$-functions, we may choose a highly ramified character $\nu_v$ such that
both $L(s, {\rm sym}^4(\pi_v)\times\pi_v\otimes\nu_v)$ and $L(s, {\rm sym}^3(\pi_v)\otimes\omega_v\nu_v)$ are both $1$. One way to see this is to use the fact that the local Langlands correspondence preserves the $L$-functions of pairs, and make use of the well known result on the Galois side. The claim now follows, thanks to (2.5).
We can deduce analogously the same statement for the $\gamma$-factors and $\varepsilon$-factors. For a representation theoretic proof of the stability of $\gamma$-factors (which yields what we want for the $\varepsilon$-factors because stability holds for the $L$-factors), see for example \cite{CoPSS},
where a general situation is treated.

Consider the global $L$-functions $L_1(s):= \, L(s, {\rm sym}^4(\pi)\times\pi)$ and
$$
L_2(s):= \, L(s, {\rm Ad}(\pi')\boxtimes\pi\otimes\omega^2)L(s, {\rm sym}^3(\pi)\otimes\omega).\leqno(2.7)
$$
They both have functional equations and analytic continuations, and moreover, thanks to (s-I3) and (2.5), they are the same local factors out side a finite set of
places containing $\Sigma_\infty$ and the ramified finite places. First choose a global character $\nu$ such that at each finite place in $S$, $\nu_v$ is sufficiently ramified so that the $\nu_v$ twists of $L_{1,v}(s)$ and $L_{2,v}(s)$ are both $1$; we can do this by the (standard) argument above. Comparing functional equations, and using the fact that the archimedean factors of $L_j(s), j=1,2$, have no poles to the right of $\Re(s)>\frac12$, we get
$$
L_{1,\infty}(s) \, = \, L_{2,\infty}(s),
$$
which furnishes (2.4).

Now pick any finite place $v$ in $S$ and choose a $\nu$ such that for every $u\in S, u \ne v$, $\nu_u$ is sufficiently ramified. Comparing the functional equations again, we get (2.2); here we use the invertibility of the $\varepsilon$-factors and the shape of the local factors, as well as the very weak Ramanujan at $v$ which separates the poles of $L_{j,v}(s)$ from its dual. After that we may, in the same way, deduce (2.3) as well. We have now proved parts (a) and (b) of the Proposition.

\medskip

Part (c) of the Proposition follows easily from the general cuspidality criterion (\cite{RaW}) for the Kim-Shahidi functorial transfer from GL$(2)\times {\rm GL}(3)$ to ${\rm GL}(6)$. In fact, already by Theorem 1.2 of \cite{Wang}, given that $\pi$ is not solvable polyhedral, which is our hypothesis, the only way Ad$(\pi')\boxtimes\pi\otimes\omega^2$ can fail to be cuspidal is for sym$^2(\pi)$ and sym$^2(\pi')$ to be twists of each other by a character. But this has been ruled out by the hypothesis, and we are done. \qed

\medskip

\section{\bf Proof of part (a) of Theorem A}

\medskip

Let $\pi$ be as in Theorem A, attached to $(\pi', \nu)$, with $\omega$ algebraic. We will show the following:

\noindent${\bf(3.1)}$
\begin{enumerate}
\item[(i)]At each place $v$, $\pi_v$ is algebraic.
\item[(ii)]The global representation $\pi_f$ is rational over a finite extension of $\Q$ (containing $\Q(\omega)$.
\end{enumerate}

\medskip

To deduce (ii) from (i), we need to show in addition that $\Q(\pi_v)$ has degree bounded
independent of the place $v$.

\medskip

\subsection{When $\pi_v$ and $\pi'_v$ are unramified}

\medskip

Let $v$ be a finite place where $\pi_v$ and $\pi'_v$ are unramified. Preserving the notations of section 1, let us recall
that the only case where we do not
yet know the algebraicity of $\pi_v$ is case (1) of Lemma 1.1, part (B). We treat this case now, making use of (sI-3).

\medskip

\noindent{\bf Lemma 3.1.1} \, \it Let $v$ be a finite place where $\pi_v$ is unramified.
Suppose we have Ad$(\pi_v)) \, \simeq$ Ad$(\pi'_v)$, which is satisfied in cases (1), (2) of Lemma 1.1, part (B). Then (sI-3) implies that one of the following identities holds:
$$
a^2=w, \, a^4=w^2, \, a^6=w^3.\leqno(3.1.2)
$$
In particular, $\pi_v$ is algebraic if $\omega_v$ is algebraic. Moreover, if $w=1$, then $a^m=1$ with $m\in\{4, 6\}$.
\rm

\medskip

\noindent{\bf Remark 3.1.2} \, The reason for also including here the case (2) of Lemma 1.1, part (B), is the following. We knew earlier that in this case, $a^2=\mu w$, with $\mu^4=1$, so that when $w=1$, $a^4=\pm 1$. But now, using Lemma 3.1.1 in addition, we rule out the (potentially troublesome) possibility $a^4=-1$ (when $w=1$), which will be important to us later.

\medskip

{\it Proof of Lemma 3.1.1}. \, Since Ad$(\pi'_v)$ is by assumption (in this Lemma) isomorphic to Ad$(\pi_v)$,
it is described by the
triple $\{{w}^{-1}a^2, 1, wa^{-2}\}$.
So Ad$(\pi'_v)\otimes\omega_v^2$ is associated to
$\{wa^2, w^2, w^3a^{-2}\}$. Thus Ad$(\pi'_v)\otimes\pi_v\otimes\omega_v^2$ is given
by the sextuple
$$
\{wa^3, w^2a, w^2a, w^3a^{-1}, w^3a^{-1}, w^4 a^{-3}\}.
$$
On the other hand, sym$^5(\pi_v)$ is attached to the sextuple
$$
\{a^5, wa^3, w^2a, w^3a^{-1}, w^4a^{-3}, w^5a^{-5}\}.
$$
Comparing these two tuples, we get
$$
\{wa^3, w^2a, w^3a^{-1}, w^4 a^{-3}\} \, = \,
\{a^5, wa^3, w^4a^{-3}, w^5a^{-5}\}.
$$
Looking at the possibilities for $w^2a$, we see that {\it either} $a^4=w^2$ (which
subsumes both the possibilities $w^2a=a^5$ and $w^2a=w^{4}a^{-3}$),
{\it or} $a^2=w$ (when $w^2a=wa^3$), {\it or} $a^6=w^3$ (when $w^2a=w^5a^{-5}$).

\qed

\medskip

Thanks to the Lemma, we see that at all the unramified places $v$, $\Q(\pi_v)$ is contained in a finite solvable extension of $\Q(\omega)$, with its degree over $\Q(\omega)$ bounded independent of $v$ (see Lemma 1.1 for cases (2) through (4). And since $\omega$ is by hypothesis (of Theorem A) algebraic, i.e., of type $A_0$,
$\Q(\omega)$ is a number field; in fact it is a CM field or totally real. Consequently, we have proved the following:

\medskip

\noindent{\bf Lemma 3.1.3} \, \it Let $\pi$ be as in Theorem A, and let $S$ be the union of the archimedean places and the finite places where $\pi_v$ is ramified. Let $\pi^S$ denote (as usual) the restricted tensor product of the local components $\pi_v$ as $v$ runs over all the places outside $S$. Then $\pi^S$ is algebraic, rational over a finite extension $K$ of $\Q(\omega)$. In fact, $K$ is contained in a compositum of cyclotomic and Kummer extensions over $\Q(\omega)$.
\rm

\medskip

\subsection{Proof of non-occurrence of $\pi_v$ of Steinberg type}

\medskip

Suppose $\pi_v$ is a special representation. Then there is a character $\lambda$ of $F_v^\ast$ such that the associated $2$-dimensional representation $\sigma_v$ of the extended Weil group $W'_{F_v}=W_{F_v}\times {\rm SL}(2, \C)$ is of the form $\lambda\otimes {st}$ where $st$ denotes the standard representation of SL$(2, \C)$. From (sI-2), it follows that $\pi'_v$ is also necessarily special, with its associated $2$-dimensional representation $\sigma'_v$ of $W'_{F_v}$ being of the form $\lambda'\otimes{st}$. Then Ad$(\pi'_v)$ has parameter Ad$(\sigma'_v)\simeq 1\otimes {\rm sym}^2(st)$, so that Ad$(\pi'_v)\otimes\pi_v\otimes\omega^2$ has the parameter
$$
\lambda^3\otimes\left({\rm sym}^3(st)\oplus st\right).\leqno(3.2.1)
$$
On the other hand, sym$^5(\pi_v)$ has the parameter
$$
\lambda^5\otimes{\rm sym}^5(st).\leqno(3.2.2)
$$
But the representations (3.2.1) and (3.2.2) cannot be isomorphic, contradicting (sI-3). Thus $\pi_v$ cannot be of Steinberg type.
\qed

\medskip

\subsection{When $\pi_v$ is unramified, but $\pi'_v$ is ramified}

\medskip

Thanks to (sI-2) and (sI-3), sym$^3(\pi'_v)$ and Ad$(\pi'_v)$ must both be unramified in this case. Suppose $\pi'_v$ is a principal series representation. So we may write
$$
\pi'_v \, \simeq \, \mu_1 \boxplus \mu_2, \, \omega'=\mu_1\mu_2,
$$
with $\mu_1^3, \mu_2^3, \mu_1^2\mu_2, \mu_1\mu_2^2, \mu_1\mu_2^{-1}$ unramified. Then $\mu_1$ and $\mu_2$ must themselves be unramified, forcing $\pi'_v$ to be in the unramified principal series.

The only remaining possibility is for $\pi'_v$ to be a supercuspidal representation with (irreducible) parameter $\sigma_v = {\rm Ind}_E^{F_v}(\chi)$, for a character $\chi$ of the multiplicative group of a quadratic extension $E$ of $F_v$. If $\theta$ denotes the non-trivial automorphism of $E/F_v$, then recall that the irreducibility of $\sigma_v$ implies that $\chi\circ\theta \ne \chi$. We have
$$
{\rm sym}^2(\sigma_v) \, \simeq \, {\rm Ind}_E^{F_v}(\chi^2)\oplus \chi_0,
$$
where $\chi_0$ is the restriction of $\chi$ to $F_v^\ast$. The determinant of $\sigma_v$, which corresponds to $\omega'_v$, is then $\chi_0\nu$, where $\nu$ is the quadratic character of $W_{F_v}$ corresponding to $E/F_v$.
Since Ad$(\pi'_v)={\rm sym}^2(\sigma_v)\otimes(\chi_0\nu)^{-1}$ is unramified, we must have $\nu$, i.e., $E/F_v$, unramified, and
$$
\chi^2\circ \theta \, = \, \chi^2,
$$
which implies
$$
\chi \, = \, \lambda(\mu\circ N_{E/F_v}), \quad {\rm with} \quad \lambda^2=1, (\lambda\circ\theta)\ne\lambda,
$$
with $\mu$ a character of $F_v^\ast$. Then $\chi_0=\lambda_0\mu^2$ and
$$
{\rm Ad}(\sigma_v) \, \simeq \, \lambda_0 \oplus\lambda_0\nu \oplus \nu.
$$
Moreover, $\lambda_0$ must be unramified as well. If $\lambda_0=1$, then ${\rm Ad}(\sigma_v)$ contains the trivial representation ${\underline{\bf 1}}$, and therefore
$$
{\rm dim}_\C \left({\rm Hom}_{W_{F_v}}({\underline{\bf 1}}, {\rm End}(\sigma_v)\right) \, = \, 2,
$$
since ${\rm End}(\sigma_v) \simeq {\rm Ad}(\sigma_v)\oplus {\underline{\bf 1}}$. This contradicts, by Schur's Lemma, the irreducibility of $\sigma_v$. So we must have $\lambda_0\ne 1$. But then, as $\lambda_0$ is an unramified quadratic character, it must coincide with $\nu$, making $\lambda_0\nu=1$. Again, we get ${\underline{\bf 1}}\subset
{\rm Ad}(\sigma_v)$, leading to a contradiction.

So we conclude that $\pi'_v$ must be unramified when $\pi_v$ is.

\medskip

\subsection{When $\pi_v$ is ramified}

\vskip 0.2in

Suppose $\pi_v$ is a ramified principal series representation attached to the characters $\mu_1, \mu_2$ of $F_v^\ast$, with $\omega_v=\mu_1\mu_2$. Then $\pi'_v$ is also necessarily a ramified principal series representation, attached to characters $\mu'_1, \mu'_2$, with $\omega'_v=\mu'_1\mu'_2$. The criteria (sI-1) and (sI-2) give conditions relating various powers of these characters and of $\nu_v$. The situation is similar to the unramified case, and we conclude algebraicity as before. Again, $\Q(\pi_v)$ is a finite extension of Kummer type over $\Q(\omega)$.

\medskip

It remains to consider the case when $\pi_v$ is supercuspidal, in which case the constraints force $\pi_v'$ to also be supercuspidal. It is well known that a supercuspidal representation with algebraic central character $\omega_v$ is algebraic, in fact rational over $\Q(\omega_v)$. Indeed, as mentioned in the Introduction, we are using the algebraic parametrization (\cite{He2}), which is equivariant for the action of Aut$(\C)$. So it suffices to verify this for the parameter $\sigma_v$ of $\pi_v$. Since we are over a local field, the image $G_v$ of $\sigma_v$ is necessarily solvable, and in particular, the image of $G_v$ in PGL$(2,\C)$ must be dihedral, tetrahedral, or octahedral. The assertion about rationality follows from the known results on the irreducible representations of coverings of $D_{2n}, A_4$ and $S_4$.

So we have now proved part (a) of Theorem A.
\qed

\medskip

\section{\bf Proof of part (b) of Theorem A}

\medskip

Here we are assuming that $\omega$ is of finite order. Then so is $\omega'$, and the arguments of sections 1 and 2 imply immediately that $\Q(\pi_f)$ is a cyclotomic field. To be precise, this is clear outside places $v$ where $\pi_v$ is not square-integrable, and the assertion in the supercuspidal case holds because it is rational over $\Q(\omega)$, which is cyclotomic.

We need to show that at any archimedean place $v$, $\pi_v$ is of Galois type. Let $\sigma_v$, resp. $\sigma'_v$, denote the $2$-dimensional representation of the Weil group $W_{F_v}$ associated to $\pi_v$, resp. $\pi'_v$. First suppose that the restriction of $\sigma_v$ to $\C^\ast$ is a one-dimensional twist of $(z/|z|)^m\oplus (z/|z|)^{-m}$ with $m>0$; this happens when either $v$ is real and $\pi_v$ is a discrete series representation (of lowest weight $m$), or $v$ is complex and $\pi_v$ is the base change of a discrete series representation of GL$(2,F_\R)$. In either case, (sI-2) implies that $\pi_v$ is also of the same form. But then we see that the restriction to $\C^\ast$ of the parameter of Ad$(\pi'_v)\otimes\pi_v$ is a one-dimensional twist of
$$
\left((z/|z|)^{2m}\oplus 1\oplus (z/|z|)^{-2m}\right)\otimes\left((z/|z|)^m\oplus(z/|z|)^{-m}\right),
$$
which is
$$
(z/|z|)^{3m}\oplus (z/|z|)^{m}\oplus((z/|z|)^m\oplus(z/|z|)^{-m}\oplus(z/|z|)^{-m}\oplus(z/|z|)^{-3m}.
$$
Since $m\ne 0$, this representation cannot possibly be a one-dimensional twist of the restriction to $\C^\ast$ of sym$^5(\sigma_v)$, which is evidently multiplicity-free. Note that we would have no contradiction if we had allowed $m=0$, which corresponds to the Galois type situation.

Thus we may assume, by the classification of unitary representations, that $\pi_v$ is a unitary character twist of a spherical representation $\sigma_v^0=\mu_1\otimes\mu_2$, where $\mu_1, \mu_2$ are unramified characters of $F_v^\ast$. Then, since $\omega_v$ has finite order by assumption, the restriction of $\sigma_v^0$ to $\C^\ast$ is necessarily of the form $\vert\cdot\vert^s\oplus\vert\cdot\vert^{-s}$. Applying (sI-2), we see that the restriction of $\sigma'_v$ to $\C^\ast$ will also need to be a unitary character twist of $\sigma_v^0{\vert_{\C^\ast}}$. Then $\left({\rm Ad}(\sigma'_v)\otimes\sigma_v\right){\vert_{\C^\ast}}$ cannot, unlike sym$^5(\sigma_v){\vert_{\C^\ast}}$, contain a unitary character twist of $\vert\cdot\vert^{5s}$, unless $s=0$. So we get a contradiction to (sI-3) if $s\ne 0$. Hence $s$ must be zero, and since $\omega_v$ has finite order, $\sigma_v$ is a finite order twist of a $\sigma_v^0$ whose restriction to $\C^\ast$ is $1 \oplus 1$. In other words, $\pi_v$ is of Galois type.

\qed

\bigskip

\section{\bf Temperedness of $\pi$}

\medskip

It suffices to prove temperedness at each place $v$. When $v$ is archimedean, we have already shown that $\pi_v$ is even of Galois type, so we may assume that $v$ is finite. Suppose $\pi_v$ is a principal series representation with parameter $\sigma_v$. If $\pi_v$ is non-tempered, i.e., violates the Ramanujan hypothesis, we must have, since $\pi_v^\vee\simeq\overline\pi_v$,
$$
\sigma_v \, \simeq \, \lambda_v\otimes\left(\vert\cdot\vert^t\oplus\vert\cdot\vert^{-s}\right), \quad {\rm with} \quad \lambda_v^2=\omega_v, \, t\in \R_+^\ast.
$$
As $\omega_v$ has finite order, $\lambda_v$ does as well. As in the archimedean spherical case, (sI-2) implies that the parameter of $\pi'_v$ is necessarily of the form
$$
\sigma'_v \, \simeq \, \lambda'_v\otimes\left(\vert\cdot\vert^t\oplus\vert\cdot\vert^{-s}\right), \quad {\rm with} \quad {\lambda'_v}^2=\omega'_v, \, t\in \R_+^\ast,
$$
for the same $t$; here we have used the fact that $\nu_v$ has finite order. It follows (since $t>0$) that
$$
{\rm Hom}\left(\vert\cdot\vert^{5t}, {\rm Ad}(\sigma'_v)\otimes\sigma_v\otimes\omega_v^2\right) \, = \, 0,
$$
while
$$
{\rm Hom}\left(\vert\cdot\vert^{5t}, {\rm sym}^5(\sigma_v)\right) \, \ne \, 0.
$$
This contradicts (sI-3) and so $t$ must be zero. In other words, $\pi_v$ must be tempered.

Finally, as is well known, if $\pi_v$ is a discrete series representation with unitary central character, it is necessarily tempered.

This proves part (c) of Theorem A.

\qed

\bigskip

\section{\bf Proof of part (d) of Theorem A}

\medskip

\noindent{\bf Proposition 6.1} \, \it Let $\pi, \pi'$ be as in Theorem A, with $\omega=\omega'=1$. Then at any finite place $v$, $\pi_v$ is either a unitary principal series or a supercuspidal representation, with $\pi_v'$ of the same type, and furthermore,
$$
\Q(\pi_v) \, \subset \, \Q[\sqrt{5}].
$$
More precisely, we have
\begin{enumerate}
\item[(a)]When $\pi_v$, $\pi_v'$ are in the unitary principal series,
\begin{enumerate}
\item[(ai)]$\Q(\pi_v) \, = \, \Q$, if Ad$(\pi_v)\simeq$ Ad$(\pi'_v)$;
\item[(aii)]$\Q(\pi_v) \, = \, \Q[\sqrt{5}]$, if Ad$(\pi_v)\not\simeq$ Ad$(\pi'_v)$.
\end{enumerate}
\item[(b)]When $\pi_v$ is supercuspidal, $\Q(\pi_v) \, = \, \Q$.
\end{enumerate}
\rm

\medskip

\noindent{\bf Corollary 6.2} \, \it Let $\pi, \pi'$ be as above (in Proposition 6.1). Then we have
$\Q(\pi_f) \, \subset \, \Q[\sqrt{5}]$, and moreover,
$$
\Q(\pi_f) = \Q[\sqrt{5}] \, \Leftrightarrow \, \pi \, \not\simeq \, \pi'.
$$
\rm

\medskip

\noindent{\bf Proposition 6.1 \, $\implies$ \, Corollary 6.2}:

\medskip

Since $\Q(\pi_f)$ is the compositum of all the $\Q(\pi_v)$ as $v$ runs over finite places (cf. \cite{C}, Prop. 3.1, for example), we have
$$
\Q(\pi_f) \, \subset \, \Q[\sqrt{5}].\leqno(6.2)
$$

First suppose Ad$(\pi_f) \, \not\simeq \, {\rm Ad}(\pi'_f)$. Then at some finite place $u$, say, Ad$(\pi_u)$ and Ad$(\pi'_u)$ need to be non-isomorphic, forcing, by (aii) of Prop. 6.1,
$$
\Q(\pi_f) = \Q(\pi_u) = \Q[\sqrt{5}].
$$

On the other hand, if Ad$(\pi_f)$ and Ad$(\pi'_f)$ are isomorphic, their local components are isomorphic as well, and (ai) of Prop. 6.1 then implies that
$$
\Q(\pi_f) = \Q.
$$

\medskip

Now we {\it claim} that for our $\pi, \pi'$,
$$
{\rm Ad}(\pi_f)\simeq {\rm Ad}(\pi'_f) \, \implies \, \pi_f \simeq \pi'_f.\leqno(6.3)
$$

Indeed, by the {\it multiplicity one theorem for SL$(2)$} (cf. \cite{Ra2}), we have
$$
\pi'_f \, \simeq \, \pi_f \otimes \nu_f,\leqno(6.4)
$$
for some idele class character $\nu$ of $F$. Since $\omega=\omega'=1$, we must have
$$
\nu^2 \, = \, 1.
$$
Anyhow, (6.4) implies
$$
{\rm sym}^3(\pi'_f) \, \simeq \, {\rm sym}^3(\pi_f)\otimes \nu_f^3,
$$
Since $\pi$ and $\pi'$ have, by (sI-2), isomorphic symmetric cubes, we must then have
$$
\nu^3 \, = \, 1,
$$
or else, sym$^3(\pi)$ will need to admit a non-trivial self-twist, which is not possible as $\pi$ is not solvable polyhedral.
So it follows that $\nu=1$, {\it proving the claim}.

This also finishes the proof of the Corollary assuming Prop. 6.1.

\qed

\bigskip

{\it Proof of Proposition 6.1}

\medskip

Since we are assuming here that $\omega=\omega'=1$, we first claim that at the places $v$ where $\pi_v$ is supercuspidal,
$\Q(\pi_v)$ is just $\Q$. Indeed, as we are using the algebraic parametrization, we can transfer the problem
to the field of rationality of the associated irreducible $2$-dimensional representation $\sigma_v$ of $W_{F_v}$. This is a known fact,
and we briefly sketch an argument. Let $\theta$ be an automorphism of $\C$. It suffices to check that $\varepsilon(s,\pi_v^\theta\otimes\chi^\theta)$
equals $\varepsilon(s,\sigma_v^\theta\otimes\chi^\theta)$, for all characters $\chi$; we are identifying, by class field theory, the characters
$\chi$ of $F_v^\ast$ with the corresponding ones of $W_{F_v}$. (We will also identify the central character $\omega_v$ with ${\rm det}(\sigma_v)$.)
If we look at the root number
$W(\sigma_v\otimes\chi) =\varepsilon(1/2,\sigma_v\otimes\chi)$, then the ratio
$W(\sigma_v^\theta\otimes\chi^\theta)/W(\omega_v^\tau\otimes\chi^\theta)$ equals
$\left(W(\sigma_v\otimes\chi)/W(\omega_v\otimes\chi)\right)^\tau$. An analogous property holds on the GL$(2)$ side, as seen by the integral representation of the local
zeta functions, using the Whittaker model, which is compatible with the $\theta$-action (\cite{Wa}. The conductors also correspond, and one gets the Galois equivariance of
$\pi_v \mapsto \sigma_v$. Next, note that we may write $\sigma_v$ as an unramified abelian twist
of an irreducible $2$-dimensional representation of Gal$(\overline F_v/F_v)$ which will have solvable, finite image.
As we saw at the end of section 3, such a $\sigma_v$ is rational over $\Q(\omega_v)$. Hence the claim.

\medskip

Moreover, we know that $\pi_v$ will not be of Steinberg type, nor non-tempered. So let us focus on the finite
places $v$ where $\pi_v$ is in the unitary principal series.

\medskip

First let $v$ be an unramified place for $\pi$. Then $w=w'=z=\delta=1$ by hypothesis, and combining Lemma 1.1 with Proposition 3.1.1, we get the following:

\medskip

\noindent{\bf Lemma 6.5} \, \it When $\pi_v, \pi'_v$ are unramified with $\omega_v=\omega'_v=1$,
one of the following happens (up to exchanging $a$ with $b=a^{-1}$ and $c$ with $d=c^{-1}$):
\begin{enumerate}
\item[(i)]$a^m=1$, for $m\in\{4,6\}$, and $c=\pm a$; Ad$(\pi_v)\simeq {\rm Ad}(\pi_v')$
\item[(ii)]$c=a^3, a=c^{-3}, a^{10}=1$; Ad$(\pi_v)\not\simeq {\rm Ad}(\pi_v')$
\end{enumerate}
\rm

\medskip

In the first case, ${\rm Ad}(\pi'_v)$ and ${\rm Ad}(\pi_v)$ are isomorphic, and $\Q(\pi_v)$ is contained in either $\Q(\zeta_4)$ or $\Q(\zeta_6)$,
where by $\zeta_n$ we mean a primitive $n$-th root of unity in $\C$. In case (ii), it is clear that $\Q(\pi_v)$ is contained in $\Q(\zeta_{10})$. On the other
hand, since $\pi_v$ is selfdual, the trace and determinant of the conjugacy class of $\pi_v$ are both real, and so $\pi_f$ is rational over the real subfield
of $\Q(\zeta_n)$ for appropriate $n$. In other words, when $\pi_v, \pi'_v$ are unramified with trivial central character,
$$
\Q(\pi_v) \, = \, \Q \quad {\rm when} \quad {\rm Ad}(\pi_v)\simeq {\rm Ad}(\pi_v') ,
$$
and
$$
\Q(\pi_v) \, = \, \Q(\sqrt{5}) \quad {\rm when} \quad {\rm Ad}(\pi_v)\not\simeq {\rm Ad}(\pi_v'),
$$

\medskip

It remains to consider the case when $\pi_v$ is a ramified, tempered principal series representation with $\omega_v=1$. Then its parameter decomposes as
$$
\sigma_v \, \simeq \, \mu \oplus \mu^{-1},
$$
for a ramified character $\mu$ of $F_v^\ast$. We have, thanks to the algebraic normalization of parameters,
$$
\Q(\pi_v) \, = \, \Q(\mu)^+,
$$
where the right hand side is the totally real subfield of the cyclotomic field $\Q(\mu)$. Again, the reason is that $\Q(\pi_v)$
is {\it a priori} contained in $\Q(\mu)$, and the selfduality of $\pi_v$ makes it rational over the real subfield.

When $\pi_v$ is such a ramified principal series representation, $\pi'_v$ is also forced to be of similar type, thanks to (sI-2), with parameter
$\sigma'_v=\mu'\oplus \mu'^{-1}$. Furthermore, we need to have
$$
\mu^3\oplus \mu\oplus \mu^{-1}\oplus \mu^{-3} \, \simeq \, {\mu'}^3\oplus \mu'\oplus {\mu'}^{-1}\oplus {\mu'}^{-3}
$$
This implies that either $\pi_v\simeq\pi'_v$ or, up to interchanging the roles of $\mu'$ and ${\mu'}^{-1}$, $\mu'=\mu^3$. In the former case, arguing as in the
unramified case above, we deduce that either $\mu^4=1$ or $\mu^6=1$. So $\Q(\pi_v)=\Q$ in this case. So we may assume that $\mu'=\mu^3$. Then (sI-3) yields
$$
\mu^5\oplus \mu^3\oplus\mu\oplus\mu^{-1}\oplus\mu^{-3}\oplus\mu^{-5} \, \simeq \, \left(\mu^6\oplus1\oplus\mu^{-6}\right)\otimes\left(\mu\oplus\mu^{-1}\right).
$$
Hence
$$
\mu^3\oplus\mu^{-3} \, \simeq \, \mu^7\oplus\mu^{-7}.
$$
In other words, we have
$$
\mu^4=1 \quad {\rm or} \quad \mu^{10}=1.
$$
Thus $\Q(\mu)$ is generated over $\Q$ by an $m$-th root of unity $\zeta_m$ with $m\in\{4, 10\}$. Consequently,
since $\Q(\zeta_4)^+=\Q$ and $\Q(\zeta_{10})^+=\Q(\sqrt{5})$,
$$
\Q(\pi_v)=\Q(\zeta_4)^+=\Q \quad {\rm or} \quad \Q(\pi_v)=\Q(\sqrt{5}).
$$

This finishes the proof of Proposition 6.1.

\qed

\medskip

Let $\tau$ be the non-trivial automorphism of $\Q(\sqrt{5})$. At each finite $v$, let $\pi_v^\tau$ denote the $\tau$-conjugate
representation of $\pi_v$, which makes sense because $\Q(\pi_v)\subset \Q(\sqrt{5})$. Similarly, let $\pi_f^\tau$ denote the $\tau$-conjugate of $\pi_f$.
Then $\pi_f^\tau$ is an admissible irreducible, generic representation because $\pi_f$ has those properties. However, it is not at all clear if $\pi_f^\tau$ is automorphic. In the present case, we can deduce the following:

\medskip

\noindent{\bf Proposition 6.6} \, \it Let $\pi$ be an s-icosahedral cusp form (with trivial central character) on GL$(2)/F$ relative to $\pi'$, also of trivial central character. Then we have
\begin{enumerate}
\item[(i)]$\pi'$ is unique, satisfying
$$
\pi_f^\tau \, \simeq \, \pi'_f.
$$
\item[(ii)]$\pi_f$ is not rational over $\Q$.
\end{enumerate} \rm

\medskip

{\it Proof}. \, We will first prove (ii) by contradiction. Suppose  $\pi_f$ is rational over $\Q$. Then by Corollary 6.2, we must have $\pi'_f\simeq\pi_f$. Hence
$$
{\rm sym}^2(\pi'_f)\otimes \pi_f \, \simeq \, {\rm sym}^3(\pi_f)\boxplus\pi_f.
$$
This means that (by (sI-3),
$$
L(s, \pi_f; {\rm sym}^5) \, = \, L(s, {\rm sym}^3(\pi_f))L(s,\pi_f)
$$
On the other hand, by Clebsch-Gordon,
$$
L(s, {\rm sym}^4(\pi_f)\times\pi_f) \, = \, L(s, \pi_f; {\rm sym}^5)L(s, {\rm sym}^3(\pi_f))
$$
So
$$
L(s, {\rm sym}^4(\pi_f)\times\pi_f) \, = \, L(s, {\rm sym}^3(\pi_f))^2L(s, \pi_f).
$$
This leads to the identity
$$
L(s, {\rm sym}^4(\pi_f)\times(\pi_f\boxtimes\pi_f)) \, = \, L(s, {\rm sym}^3(\pi_f)\times\pi_f)L(s, \pi_f\boxtimes\pi_f).
$$
The rightmost $L$-function has a pole at $s=1$ since $\pi$ is selfdual, and one knows by Shahidi that
$L(s, {\rm sym}^3(\pi_f)\times\pi_f)$ does not vanish at $s=1$. It then forces the left hand side $L$-function to have a pole at $s=1$. But it can't be, because ${\rm sym}^4(\pi)$ is, thanks to the condition (sI-1), a cusp form on GL$(5)/F$ (see \cite{KSh2}) and $\pi\boxtimes\pi$ is an automorphic form on GL$(4)/F$.
This contradiction proves that $\pi_f$ cannot be rational over $\Q$.

\medskip

To deduce (i), it suffices to prove, locally at each finite place $v$, that
$$
\pi'_v \, \simeq \, \pi_v^\tau.
$$
Let $v$ be a finite place where $\pi_v$ is in the principal series with parameter $\sigma_v=\mu\oplus\mu^{-1}$, where $\mu$ may or may not be ramified. Then $\pi'_v$ will also be of the same form, say with parameter  $\sigma'_v=\mu'\oplus{\mu'}^{-1}$. Thanks to (sI-2), we must have, up to permuting $\mu'$ and ${\mu'}^{-1}$,
$$
\mu' \, = \, \mu^m, \quad {\rm where} \quad m\in\{1,3\}.
$$
When $\mu'=\mu$, $\pi'_v\simeq\pi_v$, and as we have seen above, $\mu^m=1$ with $m\in\{4,6\}$, implying that $\Q(\pi_v)=\Q$. So let us assume that $\mu'=\mu^3$. Then (again as above) $\mu^4=1$ or $\mu^10=1$. In the former case, we get $\pi'_v\simeq\pi_v$, while in the latter case, either $\mu^2=1$, again implying $\pi'_v\simeq\pi_v$, or else,
$\mu^{10}=1$, but $\mu^2\ne 1$. Let $\tilde\tau$ be the automorphism of $M:=\Q(\zeta_{10})$ given by $\zeta_{10}\mapsto \zeta_{10}^3$. The totally real subfield of $M$ is $K=\Q(\sqrt{5})$ and $\tilde\tau$ restricts to the non-trivial automorphism $\tau$ of $K$. It follows that
$$
\pi'_v \, \simeq \, \pi_v^\tau,
$$
which holds at every principal series place $v$. In fact it holds at every finite place because at the supercuspidal places, $\Q(\pi_v)=\Q$ due to $\pi_v$ having trivial central character. Consequently, $\pi'_v\simeq \pi_v^\tau$ at each finite place $v$. Moreover, $\pi'$ is unique by the strong multiplicity one theorem.

This proves Proposition 6.6.

\qed

\medskip

We claim that $\pi'$ is s-icosahedral with $(\pi')' \simeq \pi$. Indeed, (sI-2) is a reflexive condition, and it also proves that $\pi'$ is not dihedral or tetrahedral. We get (sI-3) (for $\pi'$) by applying $\tau$ (to the (sI-3) for $\pi$) and making use of $\pi'_f$ being $\pi_f^\tau$. Furthermore, suppose sym$^4(\pi')$ is Eisensteinian. Then by \cite{KSh2}, sym$^3(\pi')$ must admit a quadratic self-twist. Then so must sym$^3(\pi)$ by (sI-2), resulting in the non-cuspidality of sym$^4(\pi)$. But this contradicts (sI-1), and the claim follows.

Finally, since $\pi'$ is s-icosahedral, it also follows, by the reasoning we employed for $\pi$, that $\pi'$ is also of Galois type at infinity.

This proves part (d) of Theorem A.

\qed

\bigskip

\noindent{\bf Remark 6.7} \, It should be noted that as $\pi$ is of Galois type at infinity with trivial central character, it cannot correspond, for $F$ totally real, to a holomorphic Hilbert modular form. So there is no reason at all, given the state of our current knowledge, to assert that its Galois conjugate $\pi^\tau$ should be modular (though it is expected). However, in our special case, it is automorphic because it happens to be isomorphic to the shadow cusp form $\pi'$.

\medskip

\medskip

\bibliographystyle{math}    
\bibliography{Icos-fibres}

\medskip

Dinakar Ramakrishnan

253-37 Caltech

Pasadena, CA 91125, USA.

dinakar@caltech.edu

\bigskip

\end{document}